%% file: authorsample.tex
\documentclass[graybox]{svmult}

\usepackage{type1cm}          
\usepackage{makeidx}          
\usepackage{graphicx}         
\usepackage{multicol}         
\usepackage[bottom]{footmisc} 
\usepackage{newtxtext}        %
\usepackage{newtxmath}        
\usepackage{indentfirst,csquotes}
\usepackage{listings}
\usepackage{mathtools}
\usepackage{adjustbox}
\usepackage{stmaryrd}
\usepackage{ulem}
\usepackage[most]{tcolorbox}
\usepackage{mdframed}
\usepackage{algorithm}
\usepackage{algpseudocode}
\usepackage{pbox}
\usepackage{xcolor}
\usepackage{hyperref}
\usepackage{etoolbox}
\usepackage{cleveref}

\input{preamble.sty}
\makeindex

\begin{document}

    \title*{Risk-averse Optimization in Random Materials: Algorithmic Advances and HPC Acceleration}
    \author{Niklas Baumgarten, Marcel Koch, David Schneiderhan, Tim Schrader, Robin Weiß}
    \institute{Niklas Baumgarten \at Heidelberg University, \email{niklas.baumgarten@uni-heidelberg.de}}

    \maketitle

    \input{src/abstract}

    \input{src/introduction}

    \input{src/risk-averse-optimization}

    \input{src/random-materials}

    \input{src/HPC-acceleration}

    \input{src/outlook-and-conclusion}

    \input{src/acknowledgment}

    \bibliographystyle{ieeetr}

    \bibliography{bibliography}

\end{document}

%% file: src/abstract.tex
\abstract{
    We summarize our advances in the algorithmic development and hardware
    utilization for risk-averse optimization problems in random materials.
    This includes risk-averse optimization using the entropic risk measure,
    as well as recently developed sampling techniques for random materials,
    that are interoperable with the optimization framework.
    Furthermore, we discuss recent progress in the efficient
    utilization of modern hybrid hardware architectures for these methods
    to solve three-dimensional partial differential equations.
}

%% file: src/introduction.tex
\section{Introduction}\label{sec:introduction}

Engineered systems that control the state of a
physical or technological process are
often subject to uncertainties arising, for example,
from incomplete or imperfect knowledge of the precise
spatial structure of the materials involved.
To account for these uncertainties,
stochastic and sample-average optimization methods,
employing sampling techniques for the underlying random materials,
have been developed that minimize a risk measure of
the discrepancy between the system state
and the desired target.
A recent survey on such problems and methods
is given in~\cite{heinkenschloss2025optimization}.

Here, we consider the entropic risk-measure
as optimization goal and impose constraints based on
partial differential equations (PDEs) in three spatial dimensions
solved with finite element methods (FEM).
The approach combines recent work on multilevel estimators
\cite{baumgarten2024fully, baumgarten2025budgeted, baumgarten2025multilevel, guth2023multilevel, nobile2025multilevel}
with convergence analysis of~\cite{guth2024parabolic, beiser2023adaptive}
resulting in a multilevel stochastic gradient
descent (MLSGD) method~\cite{baumgarten2026multilevel},
tailored to risk-averse optimization.
Risk-averse optimization seeks solutions that
perform well under uncertainty while reducing
the impact of unfavorable outcomes,
which are in the context of this paper,
random, three-dimensional materials.

Describing natural materials (e.g., soils, bones, and biological tissues)
or synthetic materials (e.g., foams, gels, composites, and concrete)
as random fields is a popular modeling approach~\cite{torquato2002random},
as it enables uncertainty quantification and provides realistic material representations.
In particular, random microstructural features such as porosity, cracks,
and inclusions can significantly influence structural integrity,
material longevity and fatigue behavior,
as well as other macroscopic properties~\cite{khristenko2020statistical,khristenko2022statistically},
and therefore must be appropriately accounted for in risk-averse applications.

While such optimization problems are highly
relevant in practice,
they also incur significant computational
costs in terms of computing time and memory
due to their three-dimensional and stochastic nature.
In the risk-averse setting, these challenges are further amplified,
since only a small number of unfavorable samples
contribute disproportionately to the optimization objective.
To mitigate these costs, we further present a new interface between
M++~\cite{baumgarten2021parallel}, the software handling the optimization
and the FEM, and Ginkgo~\cite{anzt2022ginkgo}, a linear algebra package
designed to solve large-scale linear systems on different hardware
architectures, including GPUs from different vendors.

This paper summarizes our approach to risk-averse optimization
in Section~\ref{sec:problem-statement}, introduces non-stationary and anisotropic
multiphase random material sampling based on Gaussian random fields (GRFs),
level-cut postprocessing~\cite{khristenko2020statistical, khristenko2022statistically},
and deep Gaussian processing~\cite{dunlop2018deep, latz2025deep, osborne2023convergence} in
Section~\ref{sec:random-materials}, and finally outlines implementation
details and presents initial performance evaluations of the interface between M++ and Ginkgo in 
Section~\ref{sec:hpc-acceleration}.

%% file: src/risk-averse-optimization.tex
\section{Risk-averse Optimization}\label{sec:problem-statement}

Similar to the report~\cite{baumgarten2025optimized},
we consider a PDE-constrained optimal control problem
in which the objective is formulated using the entropic risk measure.
For a random variable~$j \in L^\infty(\Omega;\mathbb{R})$
and a risk-aversion parameter~$\theta>0$,
the entropic risk is defined as
\begin{align*}
    \mathcal{R}_\theta [j]
    \coloneqq
    \frac{1}{\theta}
    \log\left(
            \int_{\Omega} \exp(\theta j(\omega)) \, \rd \PP(\omega)
    \right).
\end{align*}
The parameter~$\theta$ controls the degree of risk aversion:
as~$\theta \to 0$, the entropic risk converges to the expectation $\EE[j]$,
whereas as~$\theta \to \infty$, it approaches the essential supremum,
corresponding to a smooth worst-case optimization problem.
Hence, the entropic risk provides a continuous interpolation
between risk-neutral and worst-case optimization
making it a well suited optimization goal for SGD methods
applied to PDE-constrained optimization problems as below.

\smallskip

\textbf{Problem (Risk-averse Optimal Control of 3D Elliptic PDEs).}
Given the desired target state $\bt \in L^2(\cD)$
and a cost factor $\lambda \geq 0$,
find the optimal control $\bz \in Z$, s.t.
\begin{equation}
    \label{eq:ocp-objective}
    \min_{\bz \in Z} \, \cR_\theta \squarelr{j(\cdot, \bz)} \quad \text{with} \quad
    j(\omega, \bz) \coloneqq
    \tfrac{1}{2} \norm{\bu[\omega] - \bt}^2_{L^2(\cD)}
    + \tfrac{\lambda}{2} \norm{\bz}_{L^2(\cD)}^2 \,.
\end{equation}
The control space~$Z$ may be a subspace of~$L^2(\cD)$
equipped with box constraints, and the state
$\bu \in L^2(\Omega, H^1_0(\cD))$
is given as the solution of an elliptic PDE
on the three-dimensional domain~$\cD = (0,1)^3$.
In particular, we consider problems with a random field
$\by \in L^2(\Omega, L^{\infty}(\cD))$
representing material coefficients and homogeneous Dirichlet boundary conditions, i.e.,
\begin{equation}
    \label{eq:elliptic-pde}
    \pdeProblem{
        -\div\big(\exp(\by(\omega, \bx)) \nabla \bu(\omega, \bx) \,\big) &=& \bz(\bx) &\text{on } \,\, \Omega\times\cD \\
        \bu(\omega, \bx) &=& 0 &\text{on } \,\,  \Omega\times \partial \cD
    } \,.
\end{equation}

\begin{figure}
    \centering
    \includegraphics[width=0.55\linewidth]{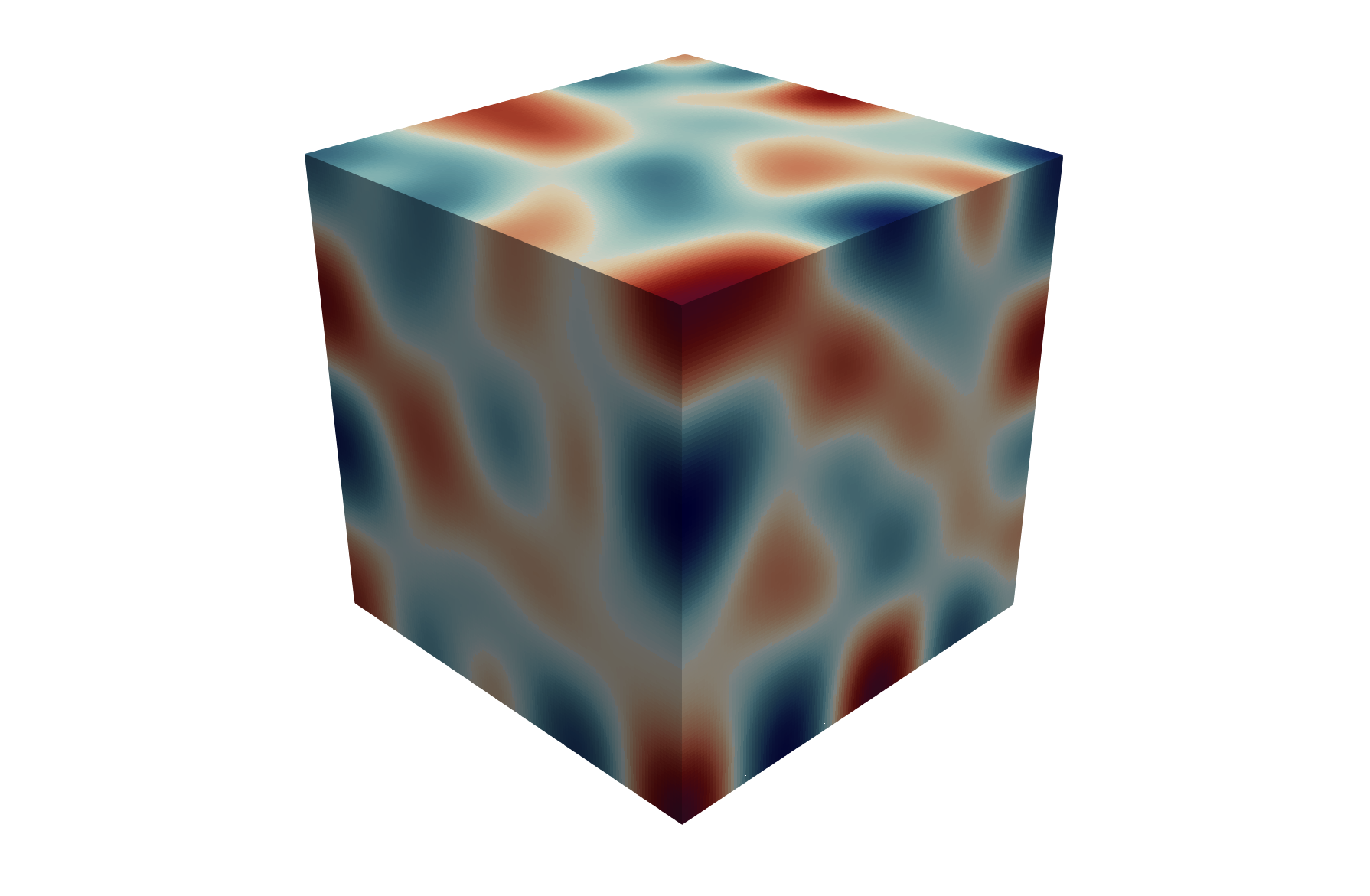}
    \hspace{-2.0cm}
    \includegraphics[width=0.55\linewidth]{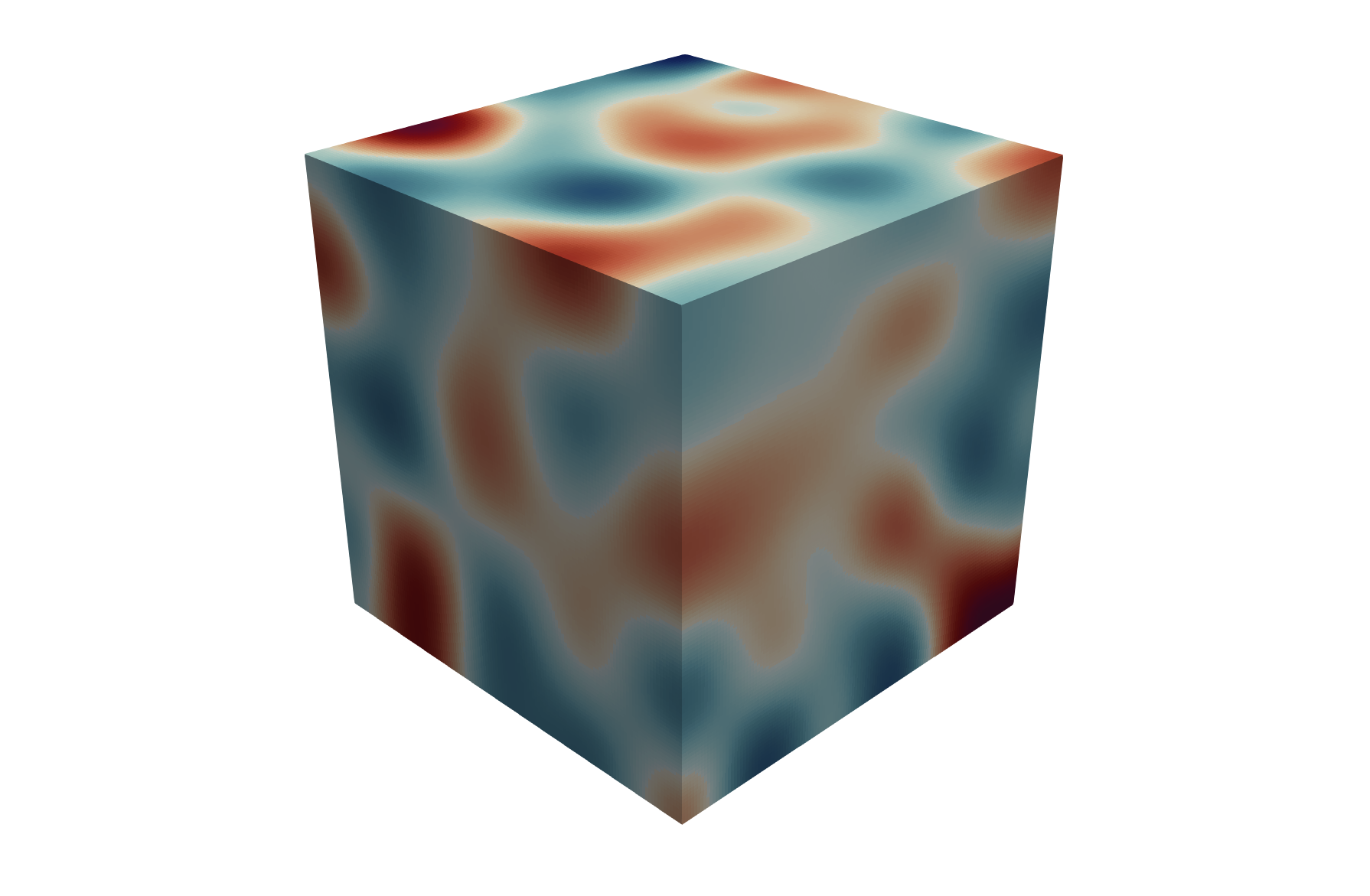}
    \caption{Samples of the stochastic fields.}
    \label{fig:samples}
\end{figure}

For our numerical tests, motivated by example~\cite[Ex.~9.37]{lord2014stochasticPDE},
we sample the realizations of the random coefficient $\by(\omega, \bx)$
with a truncated Karhunen--Loève Expansion (KLE)
\begin{equation}
    \label{eq:kle}
    \by(\omega, \bx) = \sum_{j,k,l=1}^{5} \sqrt{\lambda_{jkl}} \, \xi_{jkl}(\omega) \, \phi_{jkl}(\bx) \,,
    \quad \xi_{jkl}\sim \mathcal{U} \left(-4,4\right) \,\,\, \text{iid},
\end{equation}
using $\phi_{jkl}(\bx) = \cos(j\pi x_2) \, \cos(k\pi x_3) \, \cos(l\pi x_1)$ with $\bx = [x_1, x_2, x_3]$
and $\lambda_{jkl}=\exp(-\pi (j+k+l) \, \rho)$ with correlation length $\rho=0.15$.
Illustrations of different samples of the random field are shown in Figure~\ref{fig:samples}.
For the target state we consider $\bt(\bx) = \sin( 2 \pi x_1) \sin( 2 \pi x_2) \sin( 2 \pi x_3)$.

Problem~\eqref{eq:ocp-objective} poses a significant computational challenge,
mainly out of three reasons.
First, the underlying three-dimensional PDEs,
solved using a FEM with piecewise linear ansatz functions,
multigrid preconditioning and Krylov-space solvers,
are already computationally demanding.

Second, the optimization procedure required to approximate a
minimizer of~\eqref{eq:ocp-objective} multiplies this cost
by the number of optimization iterations and the number of sampled realizations
as it is performed through the scheme
\begin{equation}
    \label{eq:bsgd-iteration}
    \bz_\ell^{(k+1)} = \bz_\ell^{(k)} - t_k \nabla \cR_\theta[\bz_\ell^{(k)}]
\end{equation}
with $t_k$ determined by the step size rule introduced
in~\cite{koehne2024adaptivestepsizespreconditioned}
and with the gradient
\begin{equation}
    \label{eq:exact-gradient}
    \nabla \mathcal{R}_\theta[\bz_\ell^{(k)}]
    =
    \frac{\mathbb{E} \big[\!\exp\big(\theta \, \tnorm{\bu_\ell^{(k)} - \bt \,}^2_{L^2(\cD)} \big) \, \bq_\ell^{(k)} \big]}
    {\mathbb{E} \big[\!\exp\big(\theta \, \tnorm{\bu_\ell^{(k)} - \bt \,}^2_{L^2(\cD)} \big)\big]}
    + \lambda \, \bz_\ell^{(k)} \,.
\end{equation}
Here, $\bz_\ell^{(k)}$ denotes the control represented in a finite element space,
$\bu_\ell^{(k)}$ denotes the finite element approximation of the solution to~\eqref{eq:elliptic-pde} on level~$\ell$,
and $\bq_\ell^{(k)}$ analogously denotes the finite element approximation of the adjoint system
\begin{equation}
    \label{eq:elliptic-pde2}
    \pdeProblem{
        -\div\big(\exp(\by(\omega, \bx)) \nabla \bq(\omega, \bx) \,\big)
        &=& \bu(\omega, \bx) - \bt(\bx) &\text{on } \,\, \Omega\times\cD \\
        \bq(\omega, \bx) &=& 0 &\text{on } \,\, \Omega\times \partial \cD
    } \,.
\end{equation}

Third, as the risk-aversion parameter~$\theta > 0$ increases,
an increasingly large number of samples is required to accurately capture
rare but unfavorable outcomes that dominate the objective.
This can be seen in~\eqref{eq:exact-gradient},
where we proposed in~\cite{baumgarten2026multilevel}
to use multilevel Monte Carlo (MLMC) estimators
to approximate the expectation integrals.
To obtain sufficiently accurate estimates,
the number of samples required on each level increases with~$\theta$.

In conclusion, the computation of the gradient estimate~\eqref{eq:exact-gradient}
is particularly demanding, as it requires the solution of a large number
of three-dimensional PDEs.
To address these computational challenges,
especially for large values of the risk-aversion parameter~$\theta$,
we employed the HoreKa high-performance computing facility.

\smallskip

\begin{figure}
    \centering
    \includegraphics[width=1.0\linewidth]{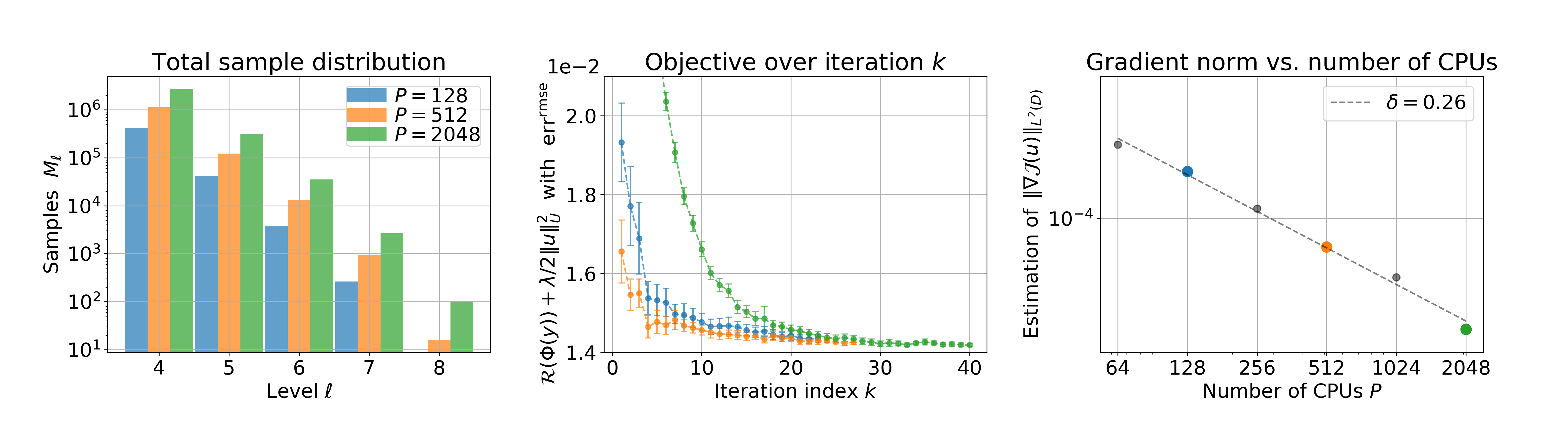}
    \caption{
        Total amount of samples used within the algorithm over MLMC levels (left),
        objective~\eqref{eq:ocp-objective} over iteration index of~\eqref{eq:bsgd-iteration} (center),
        norm of estimate to gradient estimate to~\eqref{eq:bsgd-iteration} (right).
    }
    \label{fig:numerics}
\end{figure}

\textbf{Numerical Experiments (CPU scaling).}
We outline the numerical experiments conducted in~\cite{baumgarten2026multilevel}.
A selection of the results is shown in Figure~\ref{fig:numerics},
where we fix the risk-aversion parameter to $\theta = 40$
and scale the number of CPUs from $P=64$ to $P=2048$,
marking $P=64$ in blue, $P=512$ in orange and $P=2048$ in green.

We apply the iteration scheme~\eqref{eq:bsgd-iteration},
using multilevel Monte Carlo (MLMC) to estimate the required expectations
and solving the adjoint system~\eqref{eq:elliptic-pde2}
to obtain sample-wise gradient estimates.
The overall computational budget is allocated using
techniques developed in~\cite{baumgarten2025budgeted,baumgarten2025multilevel},
resulting in approximately equal wall-clock runtimes
of four hours across all configurations.

Figure~\ref{fig:numerics} highlights three key observations.
First, the left panel shows the total number of samples,
aggregated over all optimization iterations and MLMC levels.
Both the number of samples and the number of discretization levels
used by the MLMC estimator are selected adaptively at each optimization step.
The figure clearly demonstrates that additional
computational resources are utilized effectively:
increasing the number of CPUs leads to more samples being computed and,
for $P \geq 512$, also enables the use of additional discretization levels.

As a consequence, the center panel of Figure~\ref{fig:numerics}
shows that the error estimates—represented by error bars indicating
the root-mean-square error of the gradient estimator—decrease
as more computational resources become available.
We further observe that all configurations converge to the same solution,
indicating that the parallel algorithm behaves as intended.
During development, deviations from this behavior frequently exposed implementation issues,
making this agreement an important validation of the final code.
Moreover, the increased computational resources not only permit
the computation of more samples but also accelerate batch estimation,
allowing a larger number of iterations
of~\eqref{eq:bsgd-iteration} to be performed and thereby
further reducing the optimization error.

Finally, the right panel of Figure~\ref{fig:numerics} displays
the final estimated gradient norm as a function of the number of CPUs.
We retain the color coding from the center panel and fit a
regression line through the final gradient estimates.
The results again indicate that additional computational resources
improve the quality of the solution,
since smaller gradient norms correspond to better
approximations of minimizers of~\eqref{eq:ocp-objective}.
Within the investigated range, we observe no indication of
saturation with respect to parallel scaling.
Nevertheless, the estimated convergence rate of $\delta = 0.26$
with respect to the number of CPUs is lower than the convergence
rate of $\delta = 0.31$ observed in~\cite{baumgarten2026multilevel}
with respect to wall-clock time.
We refer to~\cite{baumgarten2025multilevel,baumgarten2026multilevel,baumgarten2024fully}
for a more detailed discussion of these results.

%% file: src/random-materials.tex
\section{SPDE-Based Sampling of Random Materials}\label{sec:random-materials}

A widely used alternative to the Karhunen–Loève expansion (KLE)
as in~\eqref{eq:kle} for modeling the random input data
in~\eqref{eq:elliptic-pde} is the sampling of
Gaussian random fields (GRF) with covariance
functions from the Matérn family.
For any $\bx_1, \bx_2 \in \RR^3$,
the covariance function is given by
\begin{equation}
    \label{eq:matern-covariance}
    \Cov(\bx_1, \bx_2) = \frac{\sigma^2}{2^{\nu - 1} \Gamma(\nu)} (\kappa \br)^{\nu} \cK_{\nu}(\kappa \br), \quad
    \br = \norm{\bx_1 - \bx_2}_2, \quad \kappa = \frac{\sqrt{2\nu}}{\rho}.
\end{equation}
Here, the parameter $\nu \geq 1/2$ determines
the smoothness of the random field, $\rho > 0$
denotes the correlation length, and $\sigma > 0$
scales the variance.
Moreover, $\cK_{\nu}$ is the modified
Bessel function of the second kind,
while $\Gamma$ is the Gamma function.

GRFs with covariance
function~\eqref{eq:matern-covariance} can be
generated using the stochastic partial differential
equation (SPDE) approach introduced in~\cite{lindgren2011explicit}.
The method is based on the observation
that stationary solutions $\by(\omega, \bx)$ of the linear SPDE
\begin{equation}
    \label{eq:spde}
    \big(\kappa^2 - \Delta \big)^{\zeta} \by(\omega, \bx) = \eta \, \cW(\omega, \bx), \qquad \bx \in \RR^d
\end{equation}
are GRFs, where the randomness is driven
by Gaussian white noise $\cW$ on $\RR^d$.
The smoothness of the field is controlled
by the exponent $\zeta = \nu/2 + d/4$,
while $\eta$ is some normalization factor and
the term $(\kappa^2 - \Delta)^\zeta$ determines the
correlation structure and, in particular,
may induce isotropy or non-stationarity.
See~\cite{liu2019advances} for a survey
on further ways to sample GRFs.

In~\cite{lindgren2011explicit}, it is further proposed
to discretize~\eqref{eq:spde} using a FEM
after truncating the equation to a bounded computational domain
$\cD \subset \RR^d$.
Details on the FEM discretization of~\eqref{eq:spde},
particularly with respect to the fractional exponent $\zeta$,
the normalization $\eta$,
and the Gaussian white noise $\cW$, can be found
in~\cite{croci2018efficient, khristenko2020statistical}.
The paper~\cite{kutri2024dirichlet} introduces
Dirichlet--Neumann Averaging (DNA) as a
technique to mitigate boundary artifacts
arising from the truncation of the computational domain.
This approach is especially beneficial in memory-constrained settings,
as discussed in~\cite{baumgarten2025budgeted},
and the method of choice within this section,
even though it requires solving~\eqref{eq:spde} multiple times
with different boundary conditions.

The PDE solves required for DNA sampling are
performed using standard linear finite elements,
multigrid preconditioning, and Krylov solvers.
As discussed later in Section~\ref{subsec:performance-evaluation-and-maintenance},
these computations can also be accelerated on GPU hardware.
The resulting GRF $\by \colon \Omega \times \cD \rightarrow \RR$
(cf.~\cite[Fig.~10]{baumgarten2025optimized} for three-dimensional examples)
constitute a foundational computational output that can be further
processed to model random materials in a variety of ways.

\begin{figure}[h]
    \centering
    \includegraphics[width=1.0\linewidth]{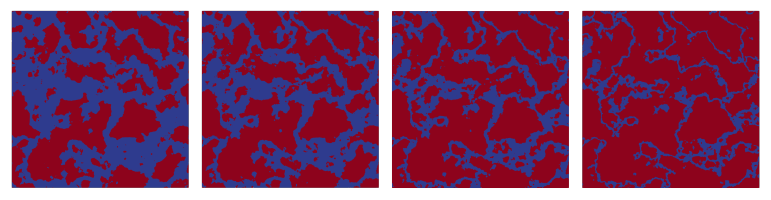}

    \vspace{-3mm}

    \caption{
        Realization of two-phase material with volume fraction $\phi_1 \in \set{0.5, \dots, 0.8}$ of red phase.
    }
    \label{fig:two-phase-2d}

    \vspace{2mm}

    \includegraphics[width=1.0\linewidth]{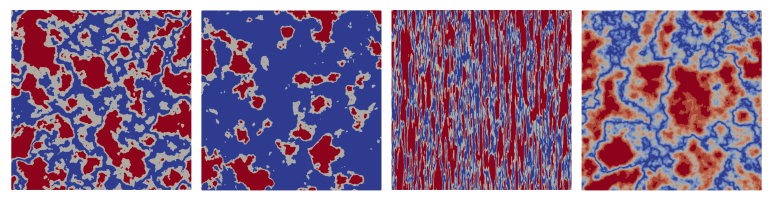}

    \vspace{-3mm}

    \caption{
        Realization of multiphase materials with volume fractions $\bphi = (\tfrac13, \tfrac13, \tfrac13)^\top$ (left),
        with volume fractions $\bphi = (0.7, 0.15, 0.15)^\top$ (center left),
        with anisotropy and volume fractions  $\bphi = (\tfrac13, \tfrac13, \tfrac13)^\top$ (center right),
        as well as ten phases $\bphi = (0.1, \dots, 0.1)^\top$ (right).
    }
    \label{fig:multi-phase-2d}

    \vspace{2mm}

    \includegraphics[width=1.0\linewidth]{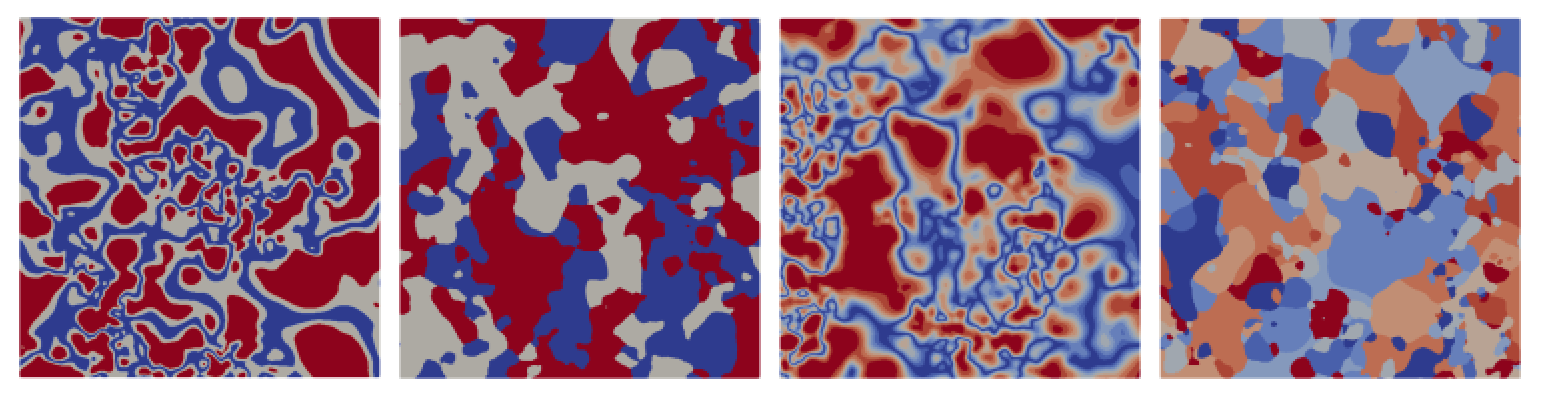}

    \vspace{-3mm}

    \caption{
        Realizations of a three-phase material (plots on the left),
        and a ten-phase material (plots on the right)
        sampled with a two-layered Gaussian process and
        ordered phases (first and third) as well as competitive
        phases (second and fourth).
    }
    \label{fig:multi-phase-deep-gaussian-2d}

\end{figure}

For example, a two-phase porous medium can be modeled
by clipping a GRF at a threshold determined by prescribed
volume fractions $\phi_1, \phi_2 \in (0,1)$ of the constituent phases such that $\phi_1 + \phi_2 = 1$.
Each point $\bx \in \cD$ is then assigned either to the
inclusion phase in $\cD_1$ or to the matrix phase in $\cD_2$,
such that $\cD_1 \cup \cD_2 = \cD$ and $\cD_1 \cap \cD_2 = \emptyset$.
See~\cite{khristenko2020statistical} for details on the construction and
Figure~\ref{fig:two-phase-2d} for an example with increasing
volume fraction of the red phase in two-spatial dimensions.

To model anisotropic covariance functions, we may extend~\eqref{eq:spde} as
\begin{align*}
    \big(\mathrm{Id} - \tfrac{1}{2 \nu} \nabla \cdot \Theta \nabla \big)^{\zeta} \by(\omega, \bx)
    = \widetilde{\eta} \, \cW(\omega, \bx), \qquad \bx \in \RR^d,
\end{align*}
where $\Theta \in \RR^{d \times d}$ is a symmetric positive-definite second-order tensor
that controls the correlation lengths
along three principal directions in a rotated reference frame,
$\mathrm{Id}$ is the identity, and $\widetilde{\eta}$
is an adjusted normalization constant~\cite{duswald2024finite}.
An extension of the amount of phases in an ordered way
by assigning volume fractions
through the vector $\bphi = (\phi_1, \dots \phi_H)$, $H \in \NN$
being the amount of total phases,
and anisotropic correlation lengths,
are illustrated in Figure~\ref{fig:multi-phase-2d}.

Furthermore, non-stationary random materials with different correlation
structure at different points in the domain can be achieved by replacing
$\kappa^2$ in~\eqref{eq:spde} with some spatially dependent function
and additional adjustments to the normalization.
Here, we even replace $\kappa^2$ with another GRF,
reassembling a deep Gaussian process, i.e,
given a fixed initial realization $\by_0(\omega, \cdot)$
we may sample a non-stationary field $\by_1 | \by_0(\omega, \cdot)$
conditioned on the initial $\by_0(\omega, \cdot)$.
For further details we refer to~\cite{dunlop2018deep,latz2025deep,osborne2023convergence}
and highlight the non-stationary correlation structures observable in
Figure~\ref{fig:multi-phase-deep-gaussian-2d}.
Here, we also illustrate competitive multiphase sampling
allowing for all phases to touch all other phases~\cite{khristenko2022statistically}
to model random imperfections and defects.

%% file: src/HPC-acceleration.tex
\section{HPC Acceleration}\label{sec:hpc-acceleration}

Finding solutions to Problem~\eqref{eq:ocp-objective}
requires solving the PDEs~\eqref{eq:elliptic-pde},~\eqref{eq:elliptic-pde2},
and if SPDE sampling for random materials is used, also~\eqref{eq:spde}.
Applying a FEM to these PDEs all result in finding
a solution $\underline{x} \in \RR^{N}$ to a sparse linear system
\begin{equation}
    \label{eq:linear-system}
    \bA \, \underline{x} = \underline{b}
    \quad \text{with} \quad
    \bA \in \RR^{N \times N} \, \text{ and } \,\, \underline{b} \in \RR^N \,.
\end{equation}
Here, $\underline{x}$ is the coefficient vector of
the finite element solution, $\bA$ is the matrix approximating
the differential operator incorporating the material coefficients
and $\underline{b}$ is the load vector capturing external forces.
Solving these systems, as there are so many of them in optimization
and uncertainty quantification (UQ) tasks,
is by far the most dominant cost.

Utilizing accelerator hardware and advanced preconditioned
Krylov solvers significantly accelerates the overall
solution process by leveraging modern HPC architectures.
As reported in~\cite{schrader_2026_19203502},
we integrated the software Ginkgo~\cite{anzt2022ginkgo}
into M++~\cite{baumgarten2021parallel} as an optional
GPU-accelerated linear algebra backend for solving~\eqref{eq:linear-system}.
M++ can now offload distributed linear algebra
operations to GPU architectures such as CUDA, HIP, and SYCL,
while taking advantage of Ginkgo’s highly configurable
solvers and preconditioners.

The motivation for this software combination is twofold:
the FEM package M++ provides deep integration of UQ and optimization methods
based on batch-wise multilevel discretization of multiple PDE samples~\cite{baumgarten2025budgeted},
which is particularly well suited for GPU architectures,
as large amounts of independent data have to be processed with a singe instruction.
At the same time, coupling both tools posed a significant challenge,
since the data layout of M++’s overlapping MPI domain decomposition
had to be reconciled with Ginkgo’s non-overlapping global indexing scheme.

\subsection{Interface design}
\label{subsec:interface-design}
\begin{figure}
    \centering
    \begin{minipage}{0.45\textwidth}
        \includegraphics[width=1.00\linewidth]{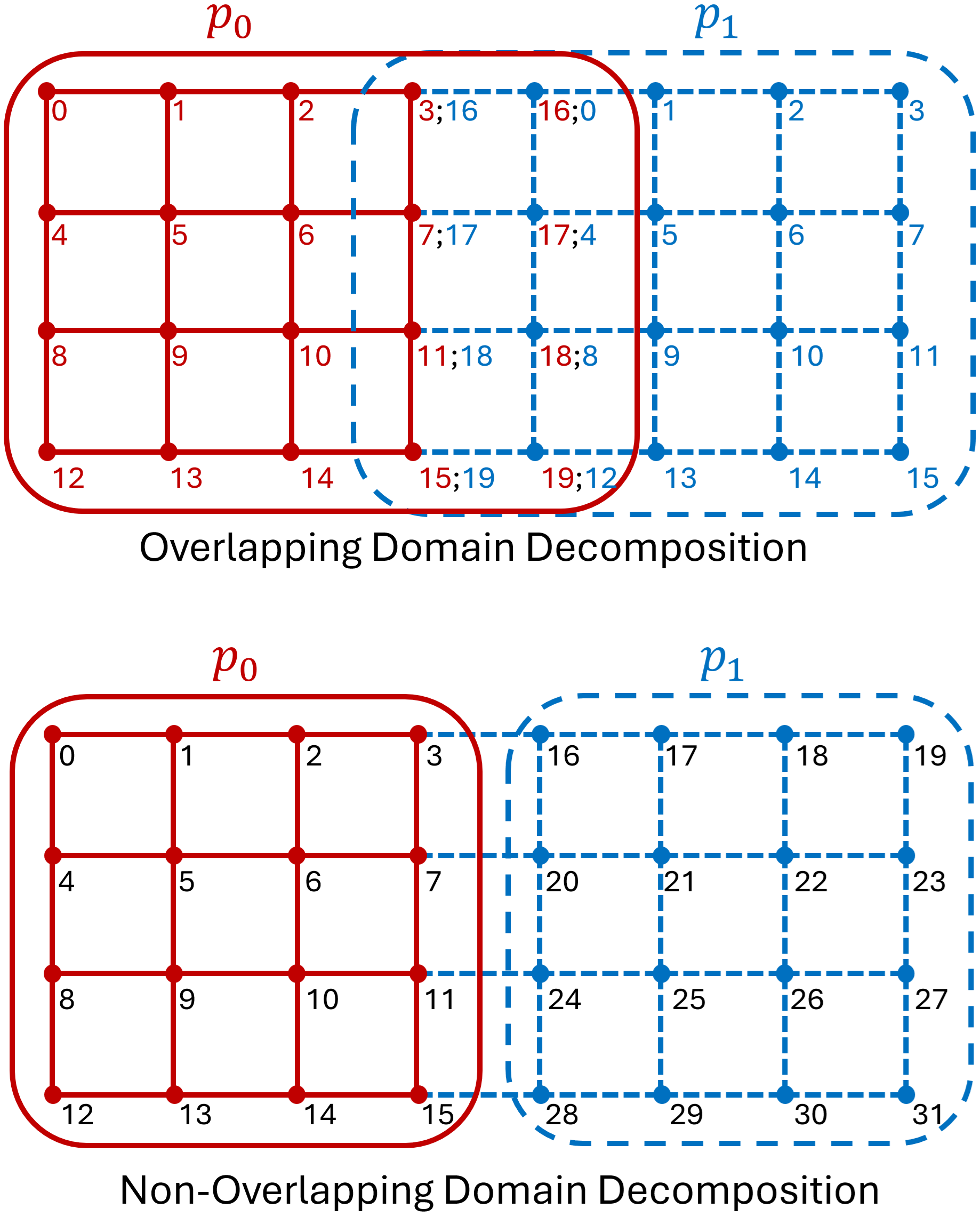}
    \end{minipage}
    \hspace{0.08\textwidth}
    \begin{minipage}{0.45\textwidth}
        \begin{align*}
            \underbrace{
                \begin{pmatrix}
                    {p_{0,0}} \\
                    {p_{0,1}} \\
                    {\vdots} \\
                    {p_{0,N}} \\
                    {p_{0,N+1}} \\
                    {\vdots} \\
                    {p_{0,N+M}} \\
                    \hline
                    {p_{1,0}} \\
                    {p_{1,1}} \\
                    {\vdots} \\
                    {p_{1,N}} \\
                    {p_{1,N+1}} \\
                    {\vdots} \\
                    {p_{1,N+M}} \\
                \end{pmatrix}
            }_{\texttt{mpp::Vector}}
            \hspace{0.5cm} \longrightarrow \hspace{-0.2cm}
            \underbrace{\begin{pmatrix}
                    {g_{0}} \\
                    {g_{1}} \\
                    {\vdots} \\
                    {g_{N}} \\
                    \ \\
                    \ \\
                    \ \\
                    \ \\
                    \hline
                    {g_{N+1}} \\
                    {g_{N+2}} \\
                    {\vdots} \\
                    {g_{2N}} \\
                    \ \\
                    \ \\
                    \ \\
                    \ \\
            \end{pmatrix}}_{\texttt{gko::distributed::Vector}}
        \end{align*}
    \end{minipage}
    \caption{
        On the mesh level (left), the overlapping domain decomposition of a 2D mesh with local indices 
        is mapped to a non-overlapping domain decomposition with global indices. 
        Color and line styles indicate rank ownership, while large rounded boxes denote the data held by a specific MPI process ($p$). 
        On the vector level (right), the locally indexed M++ vector, which includes both owned ($N$) and non-owned overlapping ($M$) points, 
        is mapped to a contiguous, globally indexed Ginkgo distributed vector.
    }
    \label{fig:mpp-gko-interface}
\end{figure}

At the heart of the new interface design is a translation layer
that maps M++’s locally indexed
and overlapping algebraic data classes, \texttt{Vector} and \texttt{Matrix},
to Ginkgo’s globally indexed dense vector representation
and compressed sparse row (CSR) matrix format.
This mapping is achieved by
separating the data within each local M++ domain into owned nodal points
(i.e., points assigned exclusively to a given rank,
including explicitly owned overlapping points)
and non-owned nodal points
(i.e., overlapping points whose data is owned by another rank).
This separation, together with the stored ownership information
for overlapping nodal points, establishes
a unified global indexing scheme for both tools
while also providing Ginkgo direct access to
the contiguously stored data owned by a rank in M++
(with copies created only when required by the target hardware memory model).
For an illustration of this translation, on how the indexing schemes refer to each other,
and on how this results in different vector representations in both tools,
we refer to Figure~\ref{fig:mpp-gko-interface}.

Beyond matching the data layouts of both tools,
the different memory models, i.e., MPI-distributed CPU memory and
multi-GPU memory, must also be homogenized.
Typically, pure MPI applications are designed to run on far more
MPI ranks $N_{\mathrm{MPI}} \in \NN$ than there are GPUs $N_{\mathrm{GPU}} \in \NN$
available in a cluster, i.e.~$N_{\mathrm{GPU}} \ll N_{\mathrm{MPI}}$.
However, a single GPU is often capable of handling the same workload as multiple CPU MPI ranks.
One way to exploit this could be to have fewer MPI ranks overall, and increase the workload per MPI rank.
This would, however, lead to increased burden on the CPUs as well.
Another approach, which is used here, is to oversubscribe the GPUs.
Each MPI rank $p_{\mathrm{MPI}} \in \NN_0$ is assigned to a GPU device identifier
$d_{\mathrm{GPU}} \in \NN_0$ according to
\begin{align*}
    d_{\mathrm{GPU}} = p_{\mathrm{MPI}} \mod N_{\mathrm{GPUs}} \,,
\end{align*}
which means that each GPU has multiple MPI ranks running on it.
The scheduling of the oversubscribed workloads on
the GPU is then handled by the underlying GPU programming model\footnote{NVIDIA provides the Multi-Process Service (MPS) for this: \url{https://docs.nvidia.com/deploy/mps/latest/index.html}.}.
The problem of over- or undersubscribing GPUs has also been discussed in~\cite{olenik2026_repartitioning},
where the authors note that two separate MPI distribution schemes,
one for the CPUs and one for the GPUs, might be even more beneficial.
In future work this could also be adapted.

Building on this reconciliation of the data layout
and memory model of the algebraic objects,
the class \texttt{GinkgoSolver},
which extends M++'s \texttt{LinearSolver},
can now directly offload algebraic computations to
Ginkgo while remaining usable for solving any PDE within M++.
This is achieved simply by selecting the executor backend
(\texttt{reference}, \texttt{omp}, \texttt{cuda}, \texttt{hip}, or \texttt{sycl})
and configuring Ginkgo either through JSON files or directly in code.

\subsection{Performance Evaluation and Maintenance}
\label{subsec:performance-evaluation-and-maintenance}
Maintaining the interface between the tools,
especially with the further advancements planned in Section~\ref{sec:outlook},
while continuously evaluating performance throughout the development process,
is achieved by integrating Ginkgo into M++’s build,
test, benchmark, and deployment pipeline~\cite{baumgarten2025continuous}.
The pipeline runs benchmarks directly on the HoreKa HPC
system\footnote{\url{https://www.nhr.kit.edu/userdocs/horeka/}},
while build and test jobs are executed in containers on
a local cluster at KIT.
As part of this setup, Ginkgo is built both in
the containers and on the HPC system, where a dedicated build
job provides it as a user space module.

\begin{figure}
    \begin{minipage}[t]{0.49\textwidth}
        \vspace{0cm}
        \centering
        \includegraphics[width=0.93\linewidth]{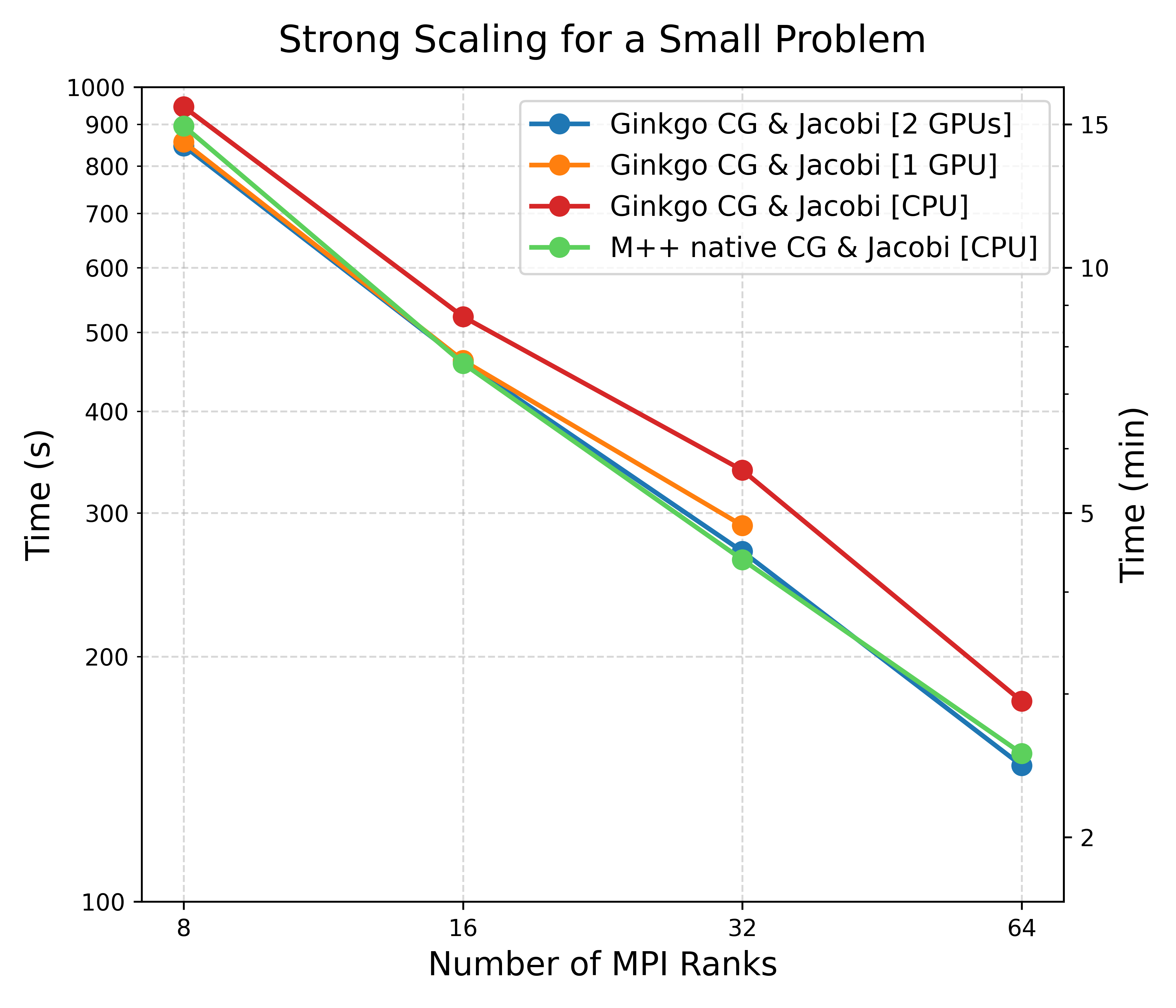}
        \caption{
            Strong scaling benchmark for a small, three dimensional problem with $5^8$ cells and 256 samples.
        }\label{fig:strong-scaling}
    \end{minipage}
    \hspace{0.2cm}
    \begin{minipage}[t]{0.49\textwidth}
        \vspace{0cm}
        \centering
        \includegraphics[width=1.0\linewidth]{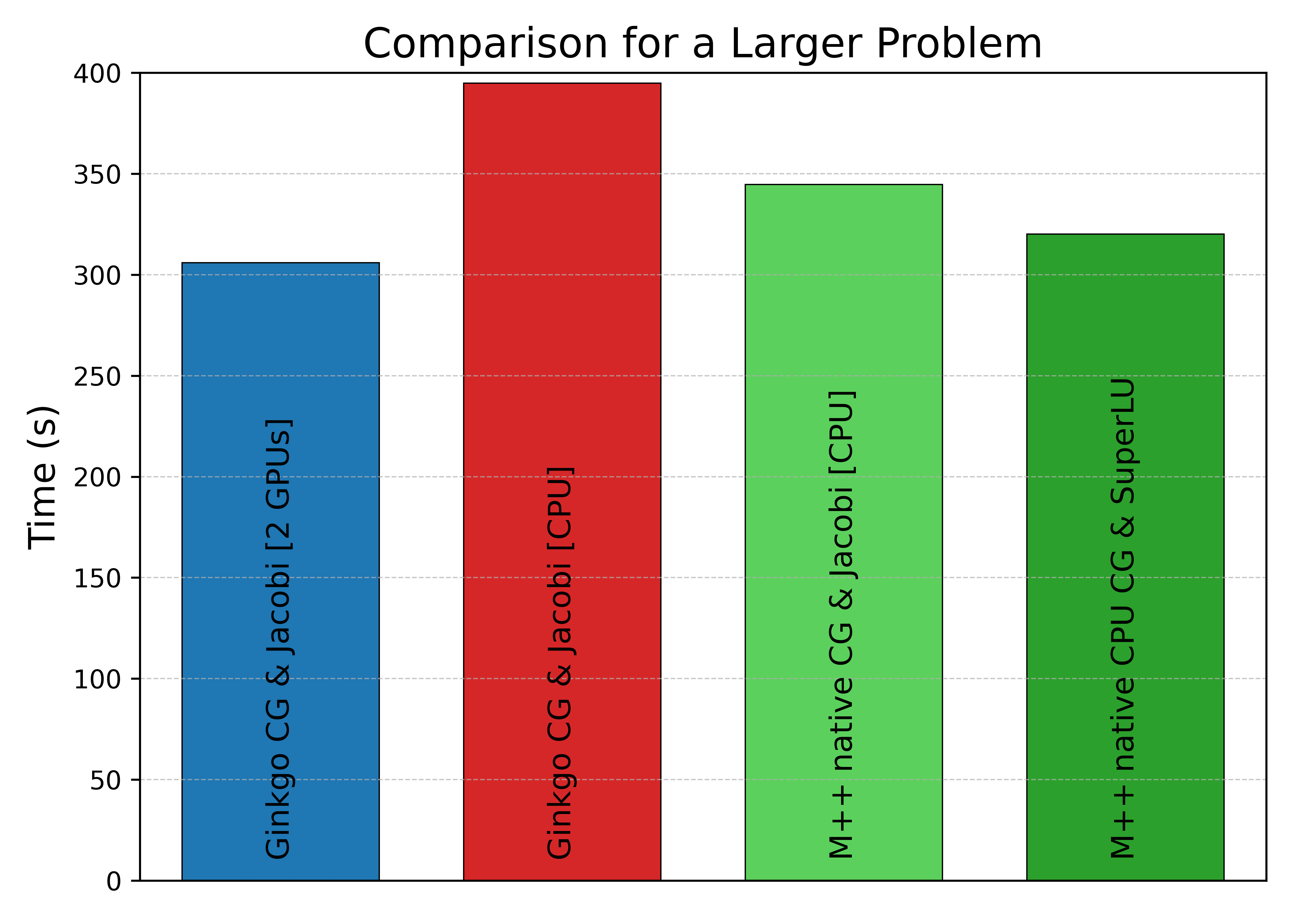}
        \vspace{0.1cm}
        \caption{
            Performance comparison for a problem with $6^8$ cells and 64 samples on 64 MPI ranks.
        }\label{fig:strong-scaling2}
    \end{minipage}
    \label{fig:figure}
\end{figure}

With this setup, maintenance is simplified, and
performance evaluations are reproducible.
As a first benchmark, cf.~Figure~\ref{fig:strong-scaling} and~\cite{schrader_2026_19203502},
we consider problem~\eqref{eq:spde},
apply SPDE based DNA sampling and discretize the arising system
with linear finite elements on $5^8$ cells.
To average out the influence of different,
random right hand sides in~\eqref{eq:spde},
we ran 256 samples resulting in 2048 linear solves
since DNA sampling in 3D requires to average over 8
different boundary conditions.

From a strong scaling experiment of this first benchmark
shown in Figure~\ref{fig:strong-scaling} running on 8 to 64 MPI ranks,
as well as from a second benchmark on $6^8$ cells,
64 samples and 64 MPI ranks shown in Figure~\ref{fig:strong-scaling2},
we can conclude four things:
i) The developed interface and targeted backends scale strongly
for these comparatively small problem and resource sizes.
ii) That one GPU performs similar to two GPUs, suggesting that
the GPU is not fully utilized.
Two GPUs are only needed to limit the oversubscription to 32-fold. 
iii) Using a Jacobi preconditioned conjugate gradient solver
in Ginkgo on CPUs via OpenMP compared to the same algorithm
implemented directly in M++ via MPI results in performance loss
around 15\% percent.
We attribute this loss, as it is a fair algorithmic and hardware comparison,
to additional steps taken in the translation layer,
which can be reduced in future developments.
iv) Despite this conversion overhead,
the GPU backend achieves a performance crossover,
resulting in an 11\% speedup over native CPU solvers
for the larger Problem in Figure~\ref{fig:strong-scaling2}.
Thus, the additional accelerator hardware overcomes the
translation loss including the additional loss from the GPU oversubscription
and still beats the pure MPI implementation, which is also using a better preconditioner,
SuperLU~\cite{li2003superlu_dist,li2023_superlu}.
All collected performance data is
available\footnote{\url{doi.org/10.5281/zenodo.19098128}}
and expanded in the future.

Preliminary tests for solving~\eqref{eq:elliptic-pde} and~\eqref{eq:elliptic-pde2}
show comparable performance between M++'s geometrical multigrid-preconditioned CG method
and Ginkgo's algebraic multigrid preconditioner running on GPUs.
Ongoing work focuses on identifying optimal configurations
for M++ when using Ginkgo as a solver, as well as improving the interface
by repartitioning the data from multiple MPI ranks onto a single rank assigned to a GPU
and reducing the amount of intermediate computation performed within the interface.

%% file: src/outlook-and-conclusion.tex
\section{Outlook and Conclusion}\label{sec:outlook}

We outlined recent advances in risk-averse optimization,
the sampling of random materials using the SPDE
approach and DNA sampling, as well as the development
of an interface connecting Ginkgo with M++ for accelerated linear algebra.
Experiments were performed on NHR@KIT’s HoreKa Blue
or Green systems\footnote{\href{https://www.nhr.kit.edu/userdocs/horeka/hardware/}{https://www.nhr.kit.edu/userdocs/horeka/hardware/}},
with a total usage
during this computing period of:
\begin{center}
        {684,717} CPU hours
    \qquad
        {133} GPU hours
\end{center}
Future work will focus particularly on expanding
the use of GPU capabilities by leveraging
the new interface described in Section~\ref{sec:hpc-acceleration}
and developing new applications involving, for example,
hyperbolic PDE systems requiring discontinuous Galerkin discretizations.

The combination and further integration of sampling-based optimization methods,
particularly those involving high-dimensional random materials,
with GPU hardware is especially promising, as computations can be performed
on independent data.
The multilevel structure, however, while reducing algorithmic complexity,
has posed a major challenge for the native M++ parallelization
and now also for the interface to Ginkgo, since problems
of vastly different sizes must be solved.
This further motivates hybrid and asynchronous approaches,
where lower levels could be treated on GPUs,
while higher levels are handled on CPUs.

%% file: src/acknowledgment.tex
\begin{acknowledgement}
    The authors gratefully acknowledge the computing time provided
    on the high-performance computer HoreKa by the National High-Performance Computing Center at KIT (NHR@KIT).
    This center is jointly supported by the Federal Ministry of Education
    and Research and the Ministry of Science, Research
    and the Arts of Baden-Württemberg,
    as part of the National High-Performance Computing (NHR)
    joint funding program.
    HoreKa is partly funded by the German Research Foundation (DFG).
\end{acknowledgement}